\documentclass{IEEEtran}
\usepackage{amssymb,latexsym,color,amsmath,pifont,epsfig,graphicx}
\usepackage{setspace,cite}
\usepackage{subfig}  
\usepackage{url}
\usepackage{rotating}
\usepackage{fancyhdr}    

\usepackage{algorithm}
\usepackage{algorithmic}

\usepackage{tikz}
\usepackage{pgfplots}
\usepackage{rotating}

\newcommand{\R}{\mathbb{R}}

\newcommand{\ds}{\displaystyle}
\newcommand{\DefinedAs}[0]{\mathrel{\mathop:}=}

\newcommand{\card}{\mbox{card}}
\newcommand{\mini}{\mathop{\mbox{minimize}}}
\newcommand{\st}{\mbox{subject to }}
\newcommand{\amin}{\mathop{\mbox{argmin}}}

\newcommand{\bff}{{\bf f}}

\newcommand{\bfh}{{\bf h}}
\newcommand{\bfu}{{\bf u}}

\newtheorem{definition}{Definition}
\newtheorem{lemma}{Lemma}
\newtheorem{proposition}{Proposition}
\newtheorem{assumption}{Assumption}
\newtheorem{remark}{Remark}

\IEEEoverridecommandlockouts                              

\begin{document}

\title{\LARGE \bf Sparse Output Feedback Synthesis via \\[0.1cm] Proximal Alternating Linearization Method}

\author{Fu Lin and Veronica Adetola
\thanks{F.\ Lin and V.\  Adetola are with Systems Department, United Technologies Research Center, 411 Silver Ln, East Hartford, CT 06118. E-mail: \{linf,adetova\}@utrc.utc.com.}
}

\maketitle

    \begin{abstract}
We consider the co-design problem of sparse output feedback and row/column-sparse output matrix. A row-sparse (resp. column-sparse) output matrix implies a small number of outputs (resp. sensor measurements). We impose row/column-cardinality constraint on the output matrix and the cardinality constraint on the output feedback gain. The resulting nonconvex, nonsmooth optimal control problem is solved by using the proximal alternating linearization method (PALM). One advantage of PALM is that the proximal operators for sparsity constraints admit closed-form expressions and are easy to implement. Furthermore, the bilinear matrix function introduced by the multiplication of the feedback gain and the output matrix lends itself well to PALM. By establishing the Lipschitz conditions of the bilinear function, we show that PALM is globally convergent and the objective value is monotonically decreasing throughout the algorithm. Numerical experiments verify the convergence results and demonstrate the effectiveness of our approach on an unstable system with 60,000 design variables.
    \end{abstract}

    
\section{Introduction}

Recent years have seen progress on the design of sparse, structured feedback controllers~\cite{linfarjov11,linfarjov12,linfarjov13,polkhlshc13,dingummicdie12,aramotkot14,wanmat14,wanlopszn14, arabahkotmot15,mat15,wanmatdoy17}. One driving force for this research direction is its wide range of applications in the control of complex systems, including power systems~\cite{dorjovchebulTPS14,wudorjovTPS16}, multi-agent systems~\cite{linfarjovALL12,zelschall13}, oscillator networks~\cite{farlinjovACC12,farlinjovTAC14sync}, and social networks~\cite{farzhalinjovCDC12}. A recent survey on the development of this research effort can be found in~\cite{jovdhiEJC16}.

A diverse set of tools for optimal sparsity control have been developed and tailored to specific design requirements. In~\cite{linfarjov11}, an augmented Lagrangian method was proposed for the structured state feedback problem. In~\cite{linfarjov12,linfarjov13}, sparse LQR state feedback controllers were obtained via the alternating direction method of multipliers. In~\cite{polkhlshc13}, an approach based on linear matrix inequality was proposed for the row/column sparse feedback problem.  In~\cite{dingummicdie12}, a convex--concave decomposition method for bilinear matrix inequality was shown effective for static output feedback problems. In~\cite{aramotkot14}, a rank constrained optimization method was developed for the sparse output feedback design. In~\cite{arabahkotmot15}, a sparse $H_2$ output feedback controller that resembles the centralized controller in frequency characteristics was proposed. In~\cite{wanmat14,mat15}, localized output feedback controllers with communication delay were developed.

In this work, we design of the output feedback and the output matrix simultaneously. The motivation for this co-design output feedback problem is two-fold. First, output feedback controllers  require fewer sensors than state feedback controllers. One may have a limited budget for the number of sensors and  is thus constrained to output feedback design. Second, it is useful to estimate the tradeoff between the number of sensors and the number of communication links for the output controllers. In practice, it is challenging to strike a good balance between the choice of sensor networks and the communication networks. Our work is a step to this direction by including the output matrix in the design process. Co-design problems of linear systems with system matrices have been considered in~\cite{chachoaza16,liuazacho16}.

The placement of sensors and actuators for feedback control has been an active research topic~\cite{lim92,gawlim96,balyou99,morpfigal15}.  In~\cite{gawlim96}, a two-part cost function was proposed for state feedback with full information and the state estimation with candidate sensors. From system integration and cost perspective, it is desired to use the least number of sensors to achieve the required performance objective~\cite{balyou99,morpfigal15}. In this context, we design simultaneously feedback sensor structure and sparse feedback gains to reduce the sensing cost and the number of communication links in distributed control.

We impose the row/column-sparsity condition on the output matrix and sparsity condition on the output feedback gain. In particular, we employ the row/column cardinality constraint in order to directly control the number of nonzero rows/columns of the output matrix. The  nonconvex, nonsmooth optimal control problem is solved by using the proximal alternating linearization method~(PALM). We establish the global convergence of PALM by proving the Lipschitz conditions of the bilinear matrix function. Furthermore, when the closed-loop performance index satisfies the Kurdyka-Lojasiewicz property, we show that PALM is guaranteed to converge to a critical point of the optimal control problem. 

The presentation is organized as follows. In Section~\ref{sec.codesign}, we formulate the co-design output feedback problem. In Section~\ref{sec.palm}, we develop the PALM algorithm and in Section~\ref{sec.conv}, we provide the convergence analysis. In Section~\ref{sec.num}, we demonstrate the convergence behavior of PALM via numerical experiments. In Section~\ref{sec.concl}, we summarize our contributions. 

\section{Co-design output feedback problem}
\label{sec.codesign}

Consider the static output feedback design
\[
\begin{array}{ll}
\dot{x}(t) \,=\, A x(t) \,+\, B_1 d(t) \,+\, B_2 u(t) \\
        y(t) \,=\, C x(t)   \\
       u(t)  \, = \, -  K y(t) 
\end{array}
\]
where $x(t) \in \R^n$ is the state, $d(t)  \in \R^q$ is the disturbance input, $u(t) \in \R^m$ is the control input, and $y(t) \in \R^p$ is the measured output.

In this work, we design both the output matrix $C \in \R^{p \times n}$ and the output feedback gain $K \in \R^{m \times p}$ simultaneously. We impose sparsity conditions on both design variables. The cardinality of the output feedback $K$ is defined as
\[
	\card(K) \, \DefinedAs \, \mbox{number of nonzero entries of $K$}.
\]
A sparser $K$ implies a smaller number of communication channels from the sensors to the actuators. We are interested in the output matrix $C$ with sparse rows or sparse columns. Because a {\em row-sparse\/} $C$ implies a small number of outputs, while a column-sparse $C$ implies a small number of sensors to measure the states. The row-cardinality of a matrix is defined as 
\begin{equation}
\nonumber
\card_{\rm row}(C) 
\,\DefinedAs\, 
 \mbox{number of nonzero rows of $C$}.
\end{equation}
Or equivalently, 
\[
\card_{\rm row}(C) \, = \, \sum_{i=1}^n \card ( \| C_i \|) ,
\] 
where $C_i$ denotes the $i$th row of $C$ and $\| \cdot \|$ denotes the Euclidean norm. Column-cardinality of $C$ is equal to the row-cardinality of its transpose, $C^T$. In what follows, we use row-sparsity without loss of generality.

The co-design problem of the sparse output feedback can be expressed as follows:
\begin{equation}
\label{eq.prob}
\begin{array}{ll}
\ds \mini_{K,C,F} & J(F) \\
\st    &    F \,=\, K C \\ 
        & \card(K) \, \leq \, s \\
        & \card_{\rm row}(C) \,\leq\, r,
\end{array}
\end{equation}
where $s$ and $r$ are prespecified positive integers. Here, $J$ is a user-specified performance index of the closed-loop system. We assume that $J$ is bounded below for all $F$. When $A - BF$ is not Hurwitz, $J$ is defined as the positive infinity. 

Problem~\eqref{eq.prob} is a nonconvex, nonsmooth optimal control problem. Because the cardinality constraints are nonconvex, nonsmooth, and the bilinear constraint $F = KC$ is nonconvex. This difficulty limits the number of solution algorithms since exiting algorithms typically require convexity or smoothness or both properties~\cite{attbol09,attbolredsou10,attbolsva13,bolsabteb14}. One may relax the cardinality constraint by using the convex surrogates such as the $\ell_1$ norm. It is noteworthy that the PALM algorithm can handle both convex and nonconvex penalty functions~\cite{bolsabteb14}.

We next put the co-design problem into a formulation suited to PALM, the proximal alternating linearization method, originally proposed for generic nonconvex, nonsmooth problems~\cite{bolsabteb14}. We begin by penalizing the difference between $F$ and $KC$ in the cost function 
\begin{equation}
\label{eq.penalty}
\begin{array}{ll}
\ds \mini_{K,C,F} &  J(F) \,+\, \dfrac{\gamma}{2} \| F \,-\, KC \|_F^2  \\
\st    
        & \card(K) \, \leq \, s , \qquad  \card_{\rm row}(C) \, \leq \, r ,
\end{array}
\end{equation}
where $\gamma$ is a sufficiently large, positive coefficient and $\|\cdot\|_F$ denotes the Frobenius norm. By introducing the indicator function
\begin{equation}
  \label{eq.f}
  f(K) \,\DefinedAs\,
  \left\{
    \begin{array}{ll}
      0, & \card(K) \,\leq\, s \\
     \infty, & \mbox{otherwise}
    \end{array}
  \right.
\end{equation}
for the cardinality constraint, and the indicator function
\begin{equation}
  \label{eq.g}
  g(C) \,\DefinedAs\,
  \left\{
    \begin{array}{ll}
      0, & \card_{\rm row}(C) \,\leq \,r \\
     \infty, & \mbox{otherwise}
    \end{array}
  \right.
\end{equation}
for the row-cardinality constraint,  problem~\eqref{eq.penalty} can be expressed as
\begin{equation}
   \label{eq.codesign}
  \mini_{K,C,F} \; \Phi \,\DefinedAs \, f(K) \,+\, g(C) \,+\, J(F) \,+\,  H(K,C,F).
\end{equation}

Note that $\Phi$ is separable with respect to $K$, $C$, and $F$ except for the coupling function
\begin{equation}
\label{eq.H}
H(K,C,F) \,\DefinedAs\, \dfrac{\gamma}{2} \| F \,-\, KC \|_F^2.
\end{equation}
It turns out that this bilinear matrix function lends itself well to PALM. In particular, the Lipschitz constants of the partial gradient of $H$ can be calculated explicitly, which facilitates the implementation of PALM and the proof of its convergence.

\section{Proximal alternating linearization method}
\label{sec.palm}

PALM falls in the class of proximal methods for nonconvex, nonsmooth optimization problems recently developed in~\cite{attbol09,attbolredsou10,attbolsva13,bolsabteb14}. It is also closely related to the alternating direction method of multipliers for convex problems. In this section, we show that the co-design problem is well suited to PALM; in particular, the proximal operators for the sparsity constraints can be computed efficiently. 

The PALM algorithm computes the minimum of the proximal functions iteratively
\begin{subequations}
\label{eq.palm}
  \begin{align}
     \label{eq.palmk}
     K^{k+1} & \;\DefinedAs\; \amin_K \left\{ f(K) \,+\, \frac{a_k}{2} \| K \,-\, X^k \|_F^2 \right\}, \\
     \label{eq.palmc}
     C^{k+1} & \;\DefinedAs\; \amin_C \left\{ g(C) \,+\, \frac{b_k}{2} \| C \,-\, Y^k \|_F^2 \right\}, \\
     \label{eq.palmf}
     F^{k+1} & \;\DefinedAs\; \amin_F \left\{ J(F) \,+\, \frac{c_k}{2} \| F \,- \, Z^k \|_F^2 \right\}.
  \end{align}
\end{subequations}
The quadratic term $\| K - X^k \|_F^2$ encourages the solution of~\eqref{eq.palmk} to be in the proximity of $X^k$, where
\[
  X^k \,=\, K^k \,-\, \frac{1}{a_k} \nabla_K H(K^k,C^k,F^k).
\]
Note that $X^k$ is a linear combination of the current iterate $K^k$ and the partial gradient of  $H$ with respect to $K$, hence the name linearization in PALM. A key requirement for the convergence of PALM is that the coefficient $a_k$ be chosen to be greater than the Lipschitz constant of $\nabla_K H$. For fixed $(C^k,F^k)$, the Lipschitz constant $L_1$ satisfies
\[
\begin{array}{l}
  \| \nabla_K H (K_1,C^k,F^k) - \nabla_K H(K_2,C^k,F^k) \|  \\[0.1cm]
  \hspace{1.5in}
 \, \leq \, L_1(C^k,F^k) \| K_1 -  K_2 \|
\end{array}
\]
for all $K_1$ and $K_2$. We set $a_k = \gamma_1 L_1$ for some $\gamma_1 > 1$. 

Similarly, the proximal points $(Y^k,Z^k)$ are linear combination of the current iterate $(C^k,F^k)$ and the partial gradients $(\nabla_C H, \nabla_F H)$,
\[
\begin{array}{l}
  Y^k \,=\, C^k \,-\, \frac{1}{b_k} \nabla_C H(K^{k+1},C^k,F^k) , \\[0.2cm]
  Z^k \,=\, F^k \,-\, \frac{1}{c_k} \nabla_F H(K^{k+1},C^{k+1},F^k).
\end{array}
\]
Let $L_2$ and $L_3$ be the Lipschitz constants of $\nabla_C H$ and $\nabla_F H$, respectively. That is, $L_2$ satisfies
\[
  \begin{array}{l}
  \| \nabla_C H (K^{k+1},C_1,F^k) - \nabla_C H(K^{k+1},C_2,F^k) \| \\[0.2cm]
    \hspace{1.5in} \leq  L_2(K^{k+1},F^k) \| C_1 -  C_2 \|
  \end{array}
\]
for all $C_1$ and $C_2$, and $L_3$ satisfies
\[
  \begin{array}{l}
  \| \nabla_F H (K^{k+1},C^{k+1},F_1) - \nabla_F H(K^{k+1},C^{k+1},F_2) \| \\[0.2cm]
    \hspace{1.5in} \leq  L_3(K^{k+1},C^{k+1}) \| F_1 -  F_2 \|
  \end{array}
\]
for all $F_1$ and $F_2$. We set 
$
b_k = \gamma_2 L_2
$
and
$
c_k = \gamma_3 L_3
$
for constants
$\gamma_2, \gamma_3 > 1$.

\subsection{Lipschitz conditions}

The Lipschitz conditions of the partial gradient of $H$ are critical for the global convergence of PALM. Furthermore, the Lipschitz conditions are necessary for the implementation of PALM because they determine the coefficients $a_k,b_k,c_k$ in the proximal operators~\eqref{eq.palm}. The Lipschitz constants for the co-design problem can be computed via a closed-form expression. This is because the partial gradient of the bilinear coupling function is linear; in particular, the partial gradients of $H$ with respect to $K$, $C$, and $F$ are given by
\[
\begin{array}{rcl}
\nabla_K H &=& \gamma (KC \,-\, F) C^T, \\
\nabla_C H &=& \gamma K^T (KC \,-\, F), \\
\nabla_F H &=& \gamma (F \,-\, KC).
\end{array}
\]
Since $\nabla_K H$, $\nabla_C H$, and $\nabla_F H$ are linear functions of $K$, $C$, and $F$, respectively, it follows that the Lipschitz constants are given by
\begin{equation}
\label{eq.lip}
L_1 \, = \, \gamma \| C C^{T}\|_F, \quad
L_2 \, = \, \gamma \| K^{T} K \|_F, \quad
L_3 \, = \, \gamma .
\end{equation}

\subsection{Explicit formulas for proximal operators}

We next show that the proximal operators~\eqref{eq.palmk} and~\eqref{eq.palmc} can be computed efficiently via explicit formulas. As a result, the implementation of PALM is particularly simple. 

The proximal operator~\eqref{eq.palmk} can be written as
\[
\begin{array}{cc}
\ds \mini_K   &  \ds \frac{a_k}{2} \| K \,-\, X^k \|_F^2  \\
\st          &  \card(K) \, \leq \, s.
\end{array}
\]
The solution is obtained by keeping the $s$ largest entries of $X^k$ in magnitude and set the remaining entries to zero.  This result is well known; e.g., see~\cite{bolsabteb14}. Let $X_s^k$ be the $s$th largest entry of $X^k$ in magnitude and let $I_s^k \in \R^{m \times n}$ be such that
\[
(I_s^k)_{ij} \,=\,
\left\{
\begin{array}{ll}
1  & \mbox{if } |X^k_{ij}| \, \geq \, X_s^k \\
0  & \mbox{otherwise}. \\
\end{array}
\right.
\]
The solution is obtained by truncating the entries whose magnitude is less than $X_s^k$ 
\begin{equation}
\label{eq.trun}
K^{k+1} \,=\, X^k \, \circ \, I_s^k,
\end{equation}
where $\circ$ denotes the entry-wise multiplication of matrices. When $\ell_1$ norm is used to promote sparsity, an efficient algorithm for the projection to the $\ell_1$ ball can be found in~\cite{ducshasin08}.

The proximal operator~\eqref{eq.palmc} can be written as
\[
\begin{array}{ll}
\ds \mini_C   &  \ds \frac{b_k}{2} \| C \,-\, Y^k \|_F^2  \\
\st          &  \card_{\rm row}(C) \, \leq \, r.
\end{array}
\]
Similar to the entry-wise truncation, the row-wise truncation amounts to keeping the $r$ largest rows of $Y^k$ in Euclidean norm and set the remaining rows to zero. Let $\delta^k$ be the $r$th largest element of $\{ \| Y_i^k \| \}_{i=1}^n $ where $Y_i^k$ denotes the $i$th row of $Y^k$. Define a binary vector $v_r^k$ of length $n$ as follows
\[
(v_r^k)_{i} \,=\,
\left\{
\begin{array}{ll}
1  & \mbox{if } \| Y^k_{i} \| \, \geq \, \delta^k \\
0  & \mbox{otherwise}. \\
\end{array}
\right.
\]
Then the row truncation of $Y^k$ can be expressed as
\begin{equation}
\label{eq.rowtrun}
      C^{k+1} \,=\, Y^k \,\circ (v_r^k {\bf 1}^T),
\end{equation} 
where ${\bf 1} \in \R^n$ is the vector of all ones. For column sparsity constraint, apply the truncation operator to the rows of $C^T$.

\subsection{Computation of proximal operator~\eqref{eq.palmf}}

One advantage of PALM is that it allows the computation of the proximal operator~\eqref{eq.palmf} to be independent of other proximal operators~\eqref{eq.palmk}-\eqref{eq.palmc}. In other words, it does not rely on a specific performance index $J$ in the minimization problem
\begin{equation}
\label{eq.minJ}
\mini_F \;\; J(F) \,+\, \frac{c_k}{2} \| F \,-\, Z^k \|_F^2.
\end{equation}
This feature of separability has been noted in the state feedback design by using the alternating direction method of multipliers~(ADMM); see~\cite{linfarjov13}.

For example, we consider the closed-loop $H_2$ norm from the disturbance $d$ to the performance output 
\[
z \,\DefinedAs \,[x^T Q^{1/2}, u^T R^{1/2}]^T,
\] 
where $Q$ and $R$ are positive definite matrices. In~\cite{linfarjov13} the Anderson-Moore method was developed for the minimization step~\eqref{eq.minJ}; see Appendix for details.

We summarize PALM in Algorithm~\ref{alg.palm}.

\begin{algorithm}[htb]
   \caption{PALM for nonconvex, nonsmooth problem~\eqref{eq.codesign}}
   \label{alg.palm}
\begin{algorithmic}
   \STATE Initialization: Start with any $(K^0,C^0,F^0)$.
   \FOR{$k=0,1,2,\ldots$ until convergence}
   \STATE // $K$-minimization step
   \STATE Compute the Lipschitz constant $L_1 = \gamma \| {C^k} C^{kT} \|_F$.
   \STATE Compute $a_k = \gamma_1 L_1(C^k)$ and the partial gradient. \\[0.1cm]
   \STATE \hspace{0.5cm} $\nabla_K H(K^k,C^k) = \gamma ( K^k {C^k} - F^k) C^{kT}$. \\[0.1cm]
   \STATE Update $X^k = K^k - \frac{1}{a_k} \nabla_K H(K^k,C^k)$.
   \STATE Perform the entry-wise truncation of $X^k$ by using~\eqref{eq.trun}.
   \STATE // $C$-minimization step
   \STATE Compute the Lipschitz constant   \\[0.1cm]
   \STATE \hspace{0.5cm} $L_2 = \gamma \| (K^{k+1})^T K^{k+1} \|_F$. \\[0.1cm]
   \STATE Compute $b_k \,=\, \gamma_2 L_2(K^{k+1})$ and the partial gradient \\[0.1cm]
   \STATE \hspace{0.5cm} $\nabla_C H(K^{k+1},C^k) \,=\, (K^{k+1})^T (K^{k+1} C^k - F^k)$. \\[0.1cm]
   \STATE Update $Y^k \,=\, C^k - \frac{1}{b_k} \nabla_C H(K^{k+1},C^k,F^k)$.
   \STATE Perform the row-wise truncation of $Y^k$  by using~\eqref{eq.rowtrun}.
   \STATE  // $F$-minimization step 
   \STATE When $J$ is the closed-loop $H_2$ norm, employ the Anderson-Moore method in Appendix~\ref{app.am}.
   \ENDFOR
\end{algorithmic}
\end{algorithm}

\section{Convergence analysis}
\label{sec.conv}

The global convergence of PALM for nonconvex, nonsmooth problems are analyzed in~\cite{attbolsva13,bolsabteb14}. In this section, we build on the results in~\cite{bolsabteb14} and show the global convergence of PALM for the co-design problem. Furthermore, the objective value is monotonically decreasing throughout the PALM algorithm. When the performance index $J$ satisfies the so-called KL property, PALM is guaranteed to converge to a critical point. The proofs can be found in Appendix.

We begin with a technical lemma on the Lipschitz conditions of $\Phi$.
\begin{lemma} \label{lem.lip}
  The objective function $\Phi$ in~\eqref{eq.codesign} satisfies the following properties:
  \begin{enumerate}
  \item \label{pro.lb}
   $\inf_{K,C,F} \Phi(K,C,F) > -\infty$, $\inf_K f(K) > -\infty$, $\inf_C g(C) > -\infty$, and $\inf_F J(F) > -\infty$.
  \item \label{pro.lip}
   The partial gradients $\nabla_K H$, $\nabla_C H$, and $\nabla_F H$ are globally Lipschitz.
\item \label{pro.lipbd}
  There exist bounded constants $q_i^-$, $q_i^+ >0$ for $i=1,2,3$ such that the Lipschitz constants in~\eqref{eq.lip} are bounded
  \begin{equation}
    \label{eq.lipbd}
    \inf_k \{ L_i^k \}  \, \geq \, q_i^-
     ~~\mbox{and}~~
    \sup_k \{ L_i^k \}  \, \leq \, q_i^+.
  \end{equation}
\item \label{pro.lipc2}
  The entire gradient $\nabla H$ is Lipschitz continuous on the bounded subsets of $\R^{m \times n} \times \R^{n \times n} \times \R^{m \times n}$.
  \end{enumerate}
\end{lemma}

Property~\ref{pro.lb}) ensures that  proximal operators in PALM are well defined and the minimization of $\Phi$ is also well defined. Property~\ref{pro.lip}) on the boundedness of the Lipschitz constants is critical for convergence. Note that the block-Lipschitz property in $K$, $C$, and $F$ is weaker than standard assumptions in proximal methods that require $\Phi$ to be globally Lipschitz in {\em joint\/} variables $(K,C,F)$; see~\cite{bolsabteb14}. Property~\ref{pro.lipbd}) guarantees that the Lipschitz constants for the partial gradients are lower and upper bounded by finite numbers. Property~\ref{pro.lipc2}) is a technical condition for controlling the distance between two consecutive steps in the sequence $(K^k,C^k,F^k)$. This is a mild condition that holds when $H$ is twice continuously differentiable as in~\eqref{eq.H}. 

\begin{assumption}	
\label{as.J}
The closed-loop performance metric $J(F) : \R^{m \times n} \to (-\infty, +\infty]$ is a proper, lower semicontinuous function.
\end{assumption}
Here $J$ is defined as the positive infinity for an unstable state feedback gain $F$. From Lemma~\ref{lem.lip} and from the convergence results established in Lemma 3.3 of~\cite{bolsabteb14} for generic PALM, it follows that the objective value $\Phi$ is monotonically decreasing in PALM. Specifically, we have the following result.
\begin{proposition}
  \label{pro.mono}
  Suppose that Assumption~\ref{as.J} holds. Let $G^k \DefinedAs (K^k,C^k,F^k)$ be a sequence generated by Algorithm~\ref{alg.palm}.  Then 
    \[
       \frac{\delta}{2} \| G^{k+1} - G^k \|^2_F \; < \; \Phi(G^k) \,-\, \Phi(G^{k+1}), \quad \forall k \geq 0
    \]
    where $\delta = \min \{ (\gamma_i - 1) q_i^- \}$ for $i=1,2,3$.
\end{proposition} 

Note that $\delta > 0$ throughout PALM iterations because $\gamma_i > 1$ for $i=1,2,3$ (see Algorithm~\ref{alg.palm}) and $q_i^- > 0$ for $i=1,2,3$ (see Lemma~\ref{lem.lip}). Thus, the convergence of the decision variable $G^k$ can be measured by the convergence of the objective value $\Phi$. The numerical experiments in Section~\ref{sec.num} verify this convergence behavior. 

Proposition~\ref{pro.mono} guarantees global convergence of PALM starting from any initial point. We next show that PALM converges to a critical point of $\Phi$ when the closed-loop performance metric satisfies the Kurdyka-Lojasiewicz~(KL) property; see Appendix~\ref{app.def} for definition. 

\begin{lemma} \label{lem.kl}
If the performance index $J$ in~\eqref{eq.codesign} satisfies the KL property, then $\Phi$ in~\eqref{eq.codesign} satisfies the KL property.
\end{lemma}

The KL property of $\Phi$ established in Lemma~\ref{lem.kl} allows us to invoke the convergence results in~\cite{bolsabteb14}.
\begin{proposition}
  \label{pro.finite}
Let $G^k = (K^k,C^k,F^k)$ be a sequence generated by Algorithm~\ref{alg.palm}. If the performance index $J$ in~\eqref{eq.codesign} satisfies the KL property, then the following results hold.
  \begin{enumerate}
    \item The sequence $\{G^k\}$ has a finite length, that is, 
    \[
        \sum_{k=1}^\infty \| G^{k+1} - G^k \|_F \; < \; \infty.
    \]
    \item The sequence $\{G^k\}$ converges to a critical point $G^* \;=\; (K^*,C^*,F^*)$ of $\Phi$.
  \end{enumerate}
\end{proposition}
Proposition~\ref{pro.finite} follows from Lemma~\ref{lem.lip} and the convergence result for generic PALM; see  Theorem 3.1 in~\cite{bolsabteb14}.

\begin{remark}[Comparison with ADMM]
The convergence analysis of ADMM typically relies on convexity assumption~\cite{boyparchueck11}. As aforementioned, no convexity assumption is required for PALM. Another noteworthy point is that ADMM is primarily used for two-block problems (i.e., two variables with a coupling constraint), while the co-design problem~\eqref{eq.codesign} is a three-block problem. It is shown in~\cite{cheheyeyua16} via a counterexample that direct extension of ADMM for multi-block convex problem may not converge. In contrast, the convergence of PALM for multi-block problems has been established in~\cite{bolsabteb14}.
\end{remark}

\begin{remark}[KL property and semi-algebraic functions]
While it may not be straightforward to establish the KL property for a given function, it is useful to show the semi-algebraic property; see Appendix~\ref{app.def} for definition. More importantly, a variety of nonsmooth functions that arise in modern applications can be shown KL via the semi-algebraic analysis, for example, all polynomial functions, indicator functions of semi-algebraic sets, finite sums and product of semi-algebraic functions, composition of semi-algebraic functions, supremum/infimum functions of semi-algebraic functions. Furthermore, several important sets are semi-algebraic, including the cone of positive semidefinite matrices, Stiefel manifolds, and matrices with constant rank. More details on the KL property and its relation to semi-algebraic functions can be found in~\cite{attbol09,attbolredsou10,attbolsva13,bolsabteb14}.
\end{remark}

\begin{remark}[Convergence rate]
Convergence rate of PALM for nonconvex, nonsmooth problems with KL property is still an on-going research topic. For semi-algebraic problems with special forms, a desingularizing technique has been developed to characterize the convergence rate. Depending on the desingularization parameters of the semi-algebraic functions, PALM converges with a finte number of steps, with a linear convergence rate, or with a sublinear rate~\cite{attbol09}. Our numerical experience suggests a linear convergence rate for the co-design output feedback problem; see Section~\ref{sec.num}. 
\end{remark}

\section{Numerical experiments}
\label{sec.num}

In this section, we illustrate the convergence property of PALM for the sparse output feedback problem. We consider a mass-spring system with $600$ design variables and an unstable system with $60,000$ design variables. For both systems PALM finds sparse solutions with prespecified sparsity levels in a few hundred steps.\footnote{It takes a few minutes on a laptop computer running Matlab 2016b with 2.4 GHz CPU and 8GB RAM.} We take the closed-loop $H_2$ norm as the performance index with $Q = I$ and $R = 10 I$. PALM is initialized with the state feedback LQR solution and the output matrix whose elements are all ones.

\subsection{Mass-spring system}

We consider the mass-spring system with $N = 10$ masses connected in series. Let $x = [\,{\bf p}^T, {\bf v}^T\,]^T$ where ${\bf p}$ and ${\bf v} \in \R^N$ denote the position and velocity of the masses, respectively. The state-space representation is given by
\[
\begin{array}{c}
A =
\left[ 
\begin{array}{cc}
O & I \\
T & O
\end{array}
\right] \in \R^{2N \times 2N},
\\[0.5cm]
B_1 = B_2 = 
\left[
\begin{array}{c}
O \\ I
\end{array}
\right] \in \R^{2N \times N}
\end{array}
\]
where $T \in \R^{N \times N}$ is a tridiagonal Toeplitz matrix with $-2$ on the main diagonal and $1$  on the first subdiagonal. 

The total number of unknown variables in $C \in \R^{20 \times 20}$ and $K \in \R^{10 \times 20}$ is 600. We set $r = N$ nonzero {\em columns\/} in $C$ and $s = 2N^2/5$ nonzero elements in $K$. In other words, we take $50\%$ column-sparsity of $C$ and $20\%$ entry-sparsity of $K$. 

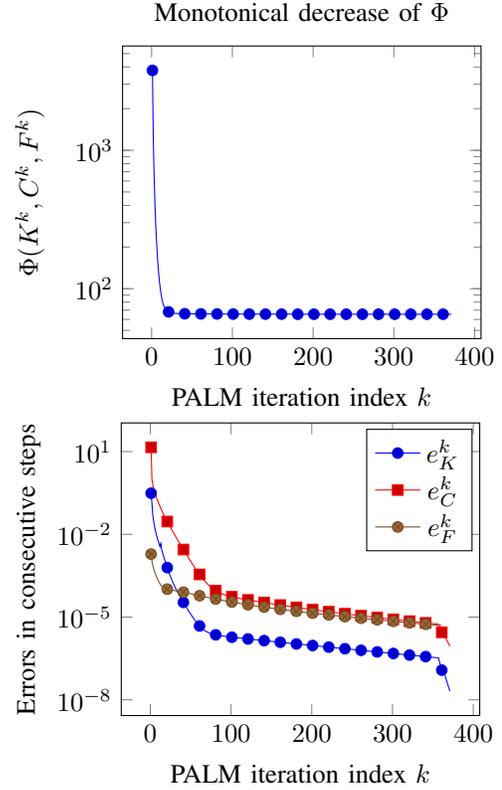
\begin{figure}
  \centering
  \begin{tikzpicture}
    \begin{semilogyaxis} [width=0.35\textwidth,
      xlabel = PALM iteration index $k$,
      mark repeat={20},
      ylabel = \mbox{$\Phi(K^k,C^k,F^k)$},
      title = Monotonical decrease of $\Phi$
      ]
      \addplot table[x=Iter,y=obj] {mass-spring-N10-colspar.txt};
    \end{semilogyaxis}
  \end{tikzpicture}
  \begin{tikzpicture}
    \begin{semilogyaxis} [width=0.35\textwidth,
      xlabel = PALM iteration index $k$,
      mark repeat={20},
      ylabel = Errors in consecutive steps
      ]
      \addplot table[x=Iter,y=diffK] {mass-spring-N10-colspar.txt};
      \addplot table[x=Iter,y=diffC] {mass-spring-N10-colspar.txt};
      \addplot table[x=Iter,y=diffF] {mass-spring-N10-colspar.txt};
      \legend{$e^k_K$,$e^k_C$,$e^k_{F}$}
    \end{semilogyaxis}
  \end{tikzpicture}
  \caption{Convergence results of PALM for the mass-spring system with column sparsity: The monotonic decreasing of $\Phi$ (top) and the convergence of the errors in two consecutive steps for the variables (bottom).}
  \label{fig.mass-spring-conv-colspar}
\end{figure}

Figure~\ref{fig.mass-spring-conv-colspar} shows the convergence results of the objective value $\Phi$ and the error of variables in consecutive  steps
\[
\begin{array}{c}
e^k_K \,=\, \| K^{k+1} - K^k \|_F,
\\[0.1cm]
e^k_C \,=\, \| C^{k+1} - C^k \|_F,
\\[0.1cm]
e^k_F \,=\, \| F^{k+1} - F^k \|_F.
\end{array}
\]
As predicted in Proposition~\ref{pro.mono}, the objective value decreases monotonically with the PALM iterations. The errors between two consecutive steps converge fast; in particular, it takes less than 300 iterations to reach $e_K^k =  4.81 \times 10^{-7}$, $e_C^k =  8.30 \times 10^{-6}$, and $e_F^k = 7.07 \times 10^{-6}$.

\begin{figure}
\centering
\begin{tabular}{cc}
\begin{tabular}{c}
\begin{sideways}
Number of inputs
\end{sideways}
\end{tabular}
&
\begin{tabular}{c}
\includegraphics[width=0.35\textwidth] {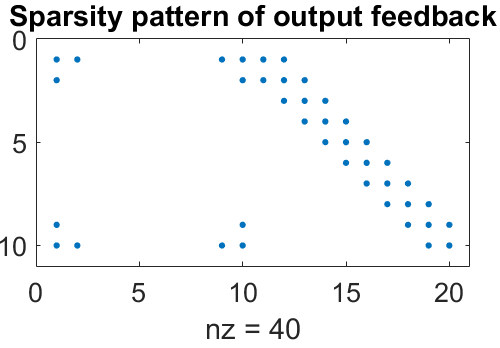}  
\end{tabular}
\\
&
\begin{tabular}{c}
Number of outputs
\end{tabular}
\end{tabular}
\\[0.1cm]
\begin{tabular}{cc}
\begin{tabular}{c}
\begin{sideways}
Number of outputs
\end{sideways}
\end{tabular}
&
\begin{tabular}{c}
\includegraphics[width=0.3\textwidth] {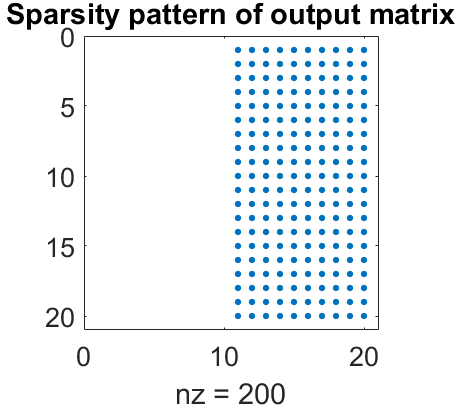}   
\end{tabular}
\\
&
\begin{tabular}{c}
Number of states
\end{tabular}
\end{tabular}
\caption{Sparsity pattern of $K$ with $20\%$ nonzero entries (top) and column-sparsity pattern of $C$ (bottom) for the mass-spring system.}
\label{fig.mass-spring-colspar}
\end{figure}

The sparsity patterns of $K$ and $C$ are shown in Fig.~\ref{fig.mass-spring-colspar}. Note that only the velocity of the masses is measured. On the other hand, the sparsity pattern of $K$ shows that the velocity of the neighboring masses is used to control the masses. The product $F = K C$ is a column sparse matrix with the same column-sparsity pattern of $C$. Therefore, one only needs $N$ sensors to measure the velocity of the masses to implement the sparse output feedback controller.

\subsection{Distributed system}

\begin{figure}
\centering
\includegraphics[width=0.25\textwidth] {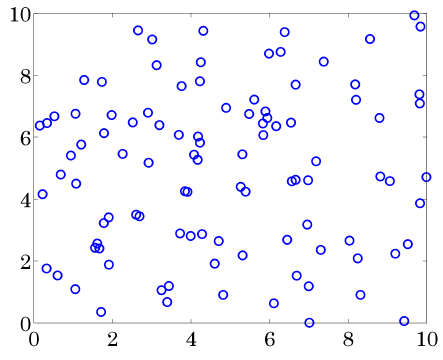}  
\caption{A network of $100$ unstable coupled systems randomly distributed in a square. The coupling strength is determined by the distance between two subsystems as in~\eqref{eq.coupled}.}
\label{fig.dist}
\end{figure}

We next consider $N = 100$ identical unstable systems in a square of $10 \times 10$ units; see Fig.~\ref{fig.dist}. The state-space representation for the $i$-th system is given by
\[
\begin{array}{rcl}
\dot{x}_i  &=& A_{ii} x_i   \,+\, \ds  \sum_{i \neq j} A_{ij} x_j   \,+\, B_i u_i  \,+\, B_i 	w_i 
\end{array}
\]
where
$
A_{ii} = 
\left[
\begin{array}{cc}
1 & 1 \\ 1 & 2
\end{array}
\right]$, 
$
A_{ij} = \alpha_{ij} I$, and
$
B_i = \left[
\begin{array}{c}
0 \\ 1
\end{array}
\right].
$
The coupling coefficient $\alpha_{ij}$ is determined by the Euclidean distance between two systems,
\begin{equation}
\label{eq.coupled}
\alpha_{ij} = 
e^{-\|p_i - p_j\|_2},
\end{equation}
where $p_i$ denotes the position of the $i$-th system. 

The total number of unknown variables in $C \in \R^{200 \times 200}$ and $K \in \R^{100 \times 200}$ is 60,000. We consider the co-design problem with $10\%$ nonzero rows in $C$ and $10\%$ nonzero entries in $K$, in other words, $r=20$ and $s=200$.

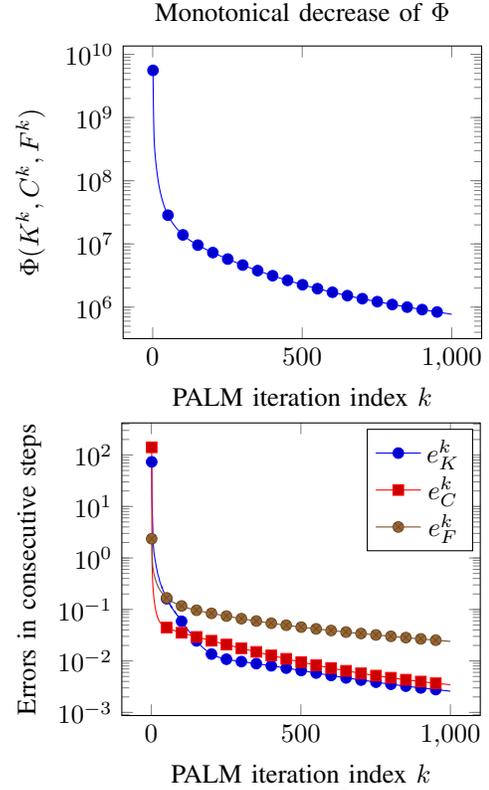
\begin{figure}
  \centering
  \begin{tikzpicture}
    \begin{semilogyaxis} [width=0.35\textwidth,
      xlabel = PALM iteration index $k$,
      mark repeat={50},
      ylabel = \mbox{$\Phi(K^k,C^k,F^k)$},
      title = Monotonical decrease of $\Phi$
      ]
      \addplot table[x=Iter,y=obj] {distsys.txt};
    \end{semilogyaxis}
  \end{tikzpicture}
  \begin{tikzpicture}
    \begin{semilogyaxis} [width=0.35\textwidth,
      xlabel = PALM iteration index $k$,
      mark repeat={50},
      ylabel = Errors in consecutive steps,
      ]
      \addplot table[x=Iter,y=diffK] {distsys.txt};
      \addplot table[x=Iter,y=diffC] {distsys.txt};
      \addplot table[x=Iter,y=diffF] {distsys.txt};
      \legend{$e^k_K$,$e^k_C$,$e^k_{F}$}
    \end{semilogyaxis}
  \end{tikzpicture}
  \caption{Convergence results of PALM for the coupled unstable system with row sparsity: The monotonic decreasing of $\Phi$ (top) and the convergence of the errors in two consecutive steps for the variables (bottom).}
  \label{fig.conv}
\end{figure}

Figure~\ref{fig.conv} shows the convergence results of the objective value $\Phi$ and the error of variables in consecutive PALM steps. As in the mass-spring example, $\Phi$ is monotonically decreasing with the PALM iterations. It takes less than 500 iterations to achieve $e_K \leq 6.6 \times 10^{-3}$,  $e_C \leq 10^{-2}$, and $e^k_F \leq 2.6 \times 10^{-2}$. 

\begin{figure}
\centering
\begin{tabular}{cc}
\begin{tabular}{c}
\begin{sideways}
Number of inputs
\end{sideways}
\end{tabular}
&
\begin{tabular}{c}
\includegraphics[width=0.4\textwidth] {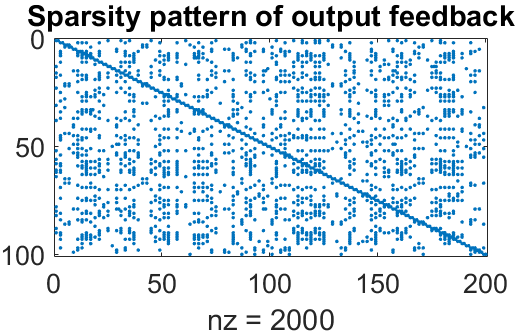}
\end{tabular}
\\
&
\begin{tabular}{c}
Number of outputs
\end{tabular}
\end{tabular}
  \\[0.25cm]
\begin{tabular}{cc}
\begin{tabular}{c}
\begin{sideways}
Number of outputs
\end{sideways}
\end{tabular}
&
\begin{tabular}{c}
\includegraphics[width=0.3\textwidth] {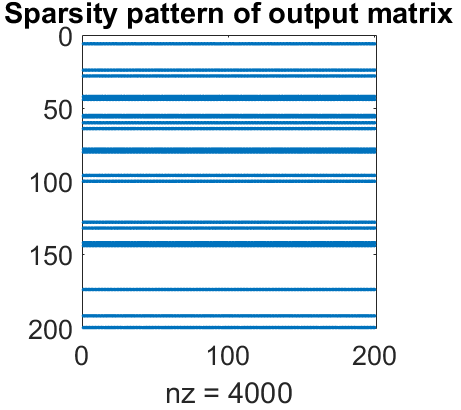}   
\end{tabular}
\\
&
\begin{tabular}{c}
Number of states
\end{tabular}
\end{tabular}
\caption{Sparsity structure of $K$ (top) and row-sparsity structure of $C$ (bottom) with $10\%$ sparsity level for the distributed system.}
\label{fig.rowspar}
\end{figure}

Figure~\ref{fig.rowspar} shows the sparsity pattern of $K$ and $C$. As required, the output matrix $C$ has exactly $r = 20$ nonzero rows ($10\%$ row-sparsity) and the output feedback gain $K$ has exactly $s = 400$ nonzero entries ($10\%$ sparsity). 

\section{Conclusions}
\label{sec.concl}

We consider the co-design problem of output feedback and output matrix simultaneously. We impose row/column-cardinality constraint to guarantee row/column sparsity on the output matrix.  We use the cardinality constraint to obtain sparse output feedback gain. The resulting nonconvex, nonsmooth problem is solved by using the PALM algorithm. We show the global convergence of PALM by establishing the Lipschitz conditions of bilinear matrix function. When the closed-loop performance index satisfies the KL property, the PALM algorithm converges to a critical point. Numerical results verify the convergence analysis and illustrate the effectiveness of our approach. 

\appendix

\subsection{Definitions of KL functions and semi-algebraic functions}
\label{app.def}

\begin{definition}[Kurdyka-Lojasiewicz property]
Let $\bff: \R^d \to (-\infty,+\infty]$ be proper and lower semicontinuous. The function $\bff$ is said to have the {\em Kurdyka-Lojasiewicz (KL)\/} property at $\bar{\bfu} \in \mbox{dom} \, \partial \bff \DefinedAs \{ \bfu \in \R^d : \partial \bff(\bfu) \neq \emptyset  \}$ if there exist $\eta \in (0, +\infty]$, a neighborhood ${\cal N}$ of $\bar{\bfu}$, and a scalar-valued function $\psi$ such that for all
\[
  \bfu \in {\cal N} \cap \{ \bff(\bar{\bfu}) < \bff(\bfu) < \bff(\bar{\bfu}) + \eta \},
\]
the following inequality holds:
\[
  \psi'(\bff(\bfu) - \bff(\bar{\bfu})) \cdot \mbox{dist}(0,\partial \bff(\bfu)) \,\geq\, 1,
\]
where $()'$ denotes the derivative function and $\mbox{dist}(x,s) \DefinedAs \inf \{ \|y-x\| : y \in {\bf s} \}$ denotes the distance from a point $x \in \R^d$ to a set ${\bf s} \subset \R^d$. A function $\bff$ is called a KL function if $\bff$ satisfies the KL property at each point of the domain of the gradient $\partial \bff$.
\end{definition}

\begin{definition}[Semi-algebraic function]
A subset ${\cal S}$ of $\R^d$ is a real {\em semi-algebraic set\/} if there exists a finite number of real polynomial functions ${\bf g}_{ij}$ and ${\bf h}_{ij} : \R^d \to \R$ such that
\[
{\cal S} \,=\, \bigcup_{j=1}^p \bigcap_{i=1}^q \{  \bfu \in \R^d: {\bf g}_{ij}(\bfu) = 0 ~\mbox{and}~ {\bf h}_{ij}(\bfu) < 0  \}.
\]
A function $\bfh:\R^d \to (-\infty,+\infty]$ is called semi-algebraic function if its graph $\{(\bfu,v) \in \R^{d+1} : \bfh(\bfu) = v \}$ is a semi-algebraic subset of $\R^{d+1}$.
\end{definition}

The connection between the KL functions and the semi-algebraic functions is provided by the following result.
\begin{proposition}[Theorem 5.1 in~\cite{bolsabteb14}]
\label{pro.semi}
A proper, lower semicontinuous, and semi-algebraic function satisfies the KL property.
\end{proposition}

\subsection{Proof of Lemma~\ref{lem.lip}}
\label{app.lip}

Property~\ref{pro.lb}) is a consequence of the coupling function $H$ in~\eqref{eq.H}, the indicator function $f$ in~\eqref{eq.f} and $g$ in~\eqref{eq.g}, and the performance metric $J$ in Assumption~\ref{as.J}. Property~\ref{pro.lip}) follows from the Lipschitz constants derived in~\eqref{eq.lip}. To show property~\ref{pro.lipbd}), $L_3(F) = \gamma$ is a constant throughout the PALM iterations. On the other hand, $L_1(C)$ in~\eqref{eq.lip} is bounded below for all $C$. Since $C^k$ is the minimizer of a feasible problem over a bounded set, it is bounded above for all $k$. Hence the entire sequence $L_1(C^k)$ satisfies the upper and lower bounds in~\eqref{eq.lipbd}. An analogous argument shows that the Lipschitz constant $L_2(K)$ satisfies~\eqref{eq.lipbd}. Finally, Property~\ref{pro.lipc2}) is a direct consequence of the twice continuous differentiability of $H$.

\subsection{Proof of Lemma~\ref{lem.kl}}
\label{app.kl}
Since KL functions are stable with respect to summation and since $J$ is assumed to be a KL function, one needs to show that $f$ and $g$ are KL functions. From Proposition~\ref{pro.semi}, we proceed to show that the indicator functions $f$ and $g$ are semi-algebraic. To this end, we use the results that the indicator function of the semi-algebraic set $\{K \,|\, \card(K) \leq s\}$ is semi-algebraic. This is because the graph of the cardinality function can be represented by a finite union of piecewise linear sets; see~\cite[Example 5.2]{bolsabteb14}. Similarly, the set of row/column sparsity matrices is also semi-algebraic. Since the indicator function of a semi-algebraic set is semi-algebraic, it follows that $f$ and $g$ are semi-algebraic functions. This completes the proof.

\subsection{Anderson-Moore method}
\label{app.am}

Let $J$ be the closed-loop $H_2$ norm from the disturbance $d$ to the performance output $z$. The necessary and sufficient conditions for the optimality of~\eqref{eq.minJ} are determined by the following coupled matrix equations~\cite{linfarjov11,linfarjov13}
\begin{subequations}
\label{eq.gradH_F}
\begin{align}
\label{eq.F}
2 (RF - B_2^T P) L \,+\, c_k (F - Z^k) \,&=\, 0 \\[0.1cm]
\label{eq.L}
(A - B_2 F) L \,+\, L (A - B_2 F)^T \,&=\, -B_1 B_1^T \\[0.1cm]
\label{eq.P}
(A - B_2 F)^T P \,+\, P (A - B_2 F) \,&=\, - (Q\,+\, F^T R F).
\end{align}
\end{subequations}
When $F$ is fixed, then~\eqref{eq.L}-\eqref{eq.P} are two Lyapunov equations in $L$ and $P$. On the other hand, when $L$ and $P$ are fixed, then~\eqref{eq.F} is a Sylvester equation in $F$. This observation motivates the Anderson-Moore method~\cite{linfarjov13}, namely, solving the Sylvester equation for $F$ and two Lyapunov equations for $(L,P)$ iteratively. The descent property of the new direction in conjunction with the Armijo line-search guarantees convergence of this approach~\cite{linfarjov12}.

The Anderson-Moore method is provided in Algorithm~\ref{alg.AM}.
\begin{algorithm}[htb]
   \caption{Anderson-Moore method}
   \label{alg.AM}
\begin{algorithmic}
   \STATE Initialization: Start with a stabilizing $F^0$
   \FOR{$l=0,1,2,\ldots$ until convergence}
   \STATE Solve~\eqref{eq.L}-\eqref{eq.P} to get the solutions $(L^l,P^l)$ .
   \STATE Solve~\eqref{eq.F} to get the solution $\bar{F}^{l}$.
   \STATE Form the direction $\Delta F^l = \bar{F}^{l} - F^l$.
   \STATE Determine stepsize $\alpha$ by using the Armijo rule.
   \STATE Update $F^{l+1} = F^l + \alpha \Delta F^l$.
   \ENDFOR
\end{algorithmic}
\end{algorithm}

\bibliographystyle{IEEEtran}
\bibliography{codesign}

\end{document}